\documentclass[12pt]{article}
\usepackage{amssymb}
\title{Le cryptosyst\`eme non-commutatif}
\author{A.Balan \footnote{balan@paris7.jussieu.fr}}
\begin{document}
\maketitle
\abstract{We show a new cryptosystem based on non-commutativ calculations
 of matrices, more specially nilpotent matrices. The cryptosystem seems powerful to restsist against usual attacks.}
\section{Introduction}
Le sujet du pr\'esent travail est de montrer l'existence d'un cryptosyst\`eme
 fond\'e sur des calculs avec des matrices nilpotentes, donc par essence
 non-commutatifs.

\bigskip

 La difficult\'e pour le d\'ecodage est celle de trouver
 la d\'ecomposition de Jordan de la matrice nilpotente. Cette derni\`ere
 est d\'etermin\'ee entre Alice et Bob par le logarithme discret modulo $p$.
 Oscar doit trouver une matrice dans un nombre consid\'erable de choix 
possibles. 
\section{Rappels d'alg\`ebre lin\'eaire}
Une matrice $X$ est nilpotente s'il existe $k$ tel que $X^k=0$ \cite{L}.

\bigskip

 Une d\'ecomposition de Jordan de la matrice permet d'obtenir une base 
 dans laquelle la matrice nilpotente de rang $n-1$ est nulle sauf sur
 la diagonale sup\'erieure o\`u elle a la valeur $1$.
\section{Le codage matriciel}
Soit une matrice nilpotente de $M_n$ de rang $n-1$ et un message $M$ en base $a$
\begin{equation}
M= \sum_i a_i a^i
\end{equation}
Le message cod\'e est la matrice suivante
\begin{equation}
\tilde M= \sum_i ln(a_i) X^i
\end{equation}
{\bf Exemple}~:
\begin{equation}
X= 
\pmatrix{0 & 1 & 0 \cr
 0 & 0 & 1 \cr
0 & 0 & 0}
\end{equation}
Pour $M=2.4 + 3.4^2$,
\begin{equation}
\tilde M=
\pmatrix{0 & 0,6932... & 1,0986... \cr
 0 & 0 & 0,6932... \cr
0 & 0 & 0}
\end{equation}
Ici le cryptage est trivial, dans le cas usuel, il faut prendre une
 matrice nilpotente autre afin de cacher les coefficients $a_i$.
\section{Le d\'ecodage du message}
Etant donn\'e $X$, on se place dans une base de Jordan de la matrice $X$,
 le message cod\'e $\tilde M$ est alors une matrice qui poss\`ede $ln(a_1)$ 
 sur la diagonale sup\'erieure, $ln(a_i)$ sur la diagonale sup\'erieure $i$.

\bigskip

Si on ne connait pas la matrice $X$, on ne peut se placer dans une telle
 base et on ne sait pas les coefficients $ln(a_k)$.
\section{D\'etermination de la matrice $X$}
On commence par calculer une matrice $A$ avec des coefficients dans
 $[0,p]$ par le syst\`eme Diffie-Hellman modulo $p$ \cite{K}\cite{Y}, en choisissant 
 la base $x$~:
\begin{equation}
{\cal A}= x^a mod(p)
\end{equation}
\begin{equation}
{\cal B}= x^b mod(p)
\end{equation}
avec $a,b$ des entiers, on d\'etermine ainsi des coefficients de la
 matrice $A$.
\begin{equation}
{\cal K}= x^{ab} mod(p) 
\end{equation}
 La matrice $X$ diff\`ere alors de la matrice standard par
 la matrice de passage suivante
\begin{equation}
e^{1/(np+\epsilon -A)}
\end{equation} 
On ajoute $\epsilon$ pour pouvoir inverser la matrice, on choisit
 ce coefficient selon les performances de l'ordinateur.

\bigskip

{\bf Exemple}~: 
\begin{equation}
{\cal A}=4^2= 2 mod(7),
\end{equation}
\begin{equation}
{\cal B}=4^4=4 mod(7),
\end{equation}
\begin{equation}
{\cal K}=4^8= 2 mod(7)
\end{equation}
\section{R\'esistance du cryptosyst\`eme}
Le cryptosyst\`eme, pour un d\'ecodage simple demanderait
\begin{equation}
e^{n^2(1/\epsilon+ln(a))}
\end{equation}
op\'erations. La difficult\'e est donc grande si on emploie une m\'ethode
 simple pour le d\'ecryptage; pour calculer $X$, il faudrait r\'esoudre
 le logarithme discret $n^2$ fois, pour chaque coefficient de la matrice
 $A$. 


\begin{thebibliography}{ccc}
\bibitem{B}
J.Buchmann, {\it Introduction to Cryptography}, Springer-Verlag, 2000.
\bibitem{}
M.A.Cherepnev, {\it On the connection between the dicrete logarithms and
 the Diffie-Hellman problem}, DiskretnajaMat.,6, 1996, pp341-349.
\bibitem{K}
N.Koblitz, {\it A course in Number Theory and Cryptography}, 2ed, Springer-Verlag, 1993.
\bibitem{L}
S.Lang, {\it Algebra}, 2ed, Addison-Wesley, 1984.
\bibitem{M}
K.S.McCurley, {\it The Dicrete Logarithm Problem}, Cryptology and
 Computional Number Theory, Proceedings of Symposia in Applied Mathematics
 42, AMS, 1990, pp49-74.
\bibitem{S}
I.Shparlinski, {\it Number Theoretic Methods in Cryptography},
 Birkhaeuser, 1999.
\bibitem{Y}
S.Y.Yan, {\it Number Theory for Computing}, Springer-Verlag, 2000.
\end{thebibliography}
\end{document}